\documentclass[twoside,11pt,reqno]{amsart}
\usepackage{amsmath,amssymb,amscd}

\makeatletter

\@settopoint\oddsidemargin
\@settopoint\marginparwidth
\@settopoint\evensidemargin
\@settopoint\topmargin
\newtheorem{Proposition}{Proposition}[section]
\newtheorem{Lemma}[Proposition]{Lemma}
\newtheorem{Theorem}[Proposition]{Theorem}
\newtheorem{Corollary}[Proposition]{Corollary}

\@addtoreset{equation}{section}

\def\F{{\mathrm F}}
\def\C{{\mathbb C}}

\def\Z{{\mathbb Z}}

\def\coeff{{\operatorname{coeff}}}
\def\sgn{{\operatorname{sgn}}}
\def\ad{{\operatorname{ad}}}
\def\Ad{{\operatorname{Ad}}}
\def\cdet{\operatorname{cdet}}

\def\row{\operatorname{row}}
\def\col{\operatorname{col}}
\def\gr{\operatorname{gr}}

\newdimen\hoogte    \hoogte=12pt    
\newdimen\breedte   \breedte=12pt   
\newdimen\dikte     \dikte=0.5pt    

\newenvironment{young}{\begingroup
       \def\vr{\vrule height0.8\hoogte width\dikte depth 0.2\hoogte}
       \def\fbox##1{\vbox{\offinterlineskip
                    \hrule height\dikte
                    \hbox to \breedte{\vr\hfill##1\hfill\vr}
                    \hrule height\dikte}}
       \vbox\bgroup \offinterlineskip \tabskip=-\dikte \lineskip=-\dikte
            \halign\bgroup &\fbox{##\unskip}\unskip  \crcr }
       {\egroup\egroup\endgroup}
\def\diagram#1{\relax\ifmmode\vcenter{\,\begin{young}#1\end{young}\,}\else%
              $\vcenter{\,\begin{young}#1\end{young}\,}$\fi}

\def\mf {\mathfrak}

\begin{document}
\title[Centralizers of nilpotent matrices]{Elementary invariants for centralizers
of nilpotent matrices}
\author{Jonathan Brown and Jonathan Brundan}
\address{Department of Mathematics, University of Oregon, Eugene, OR 97403.}
\email{jbrown8@uoregon.edu, brundan@uoregon.edu}
\thanks{Research supported in part by NSF grant no. DMS-0139019.}

\begin{abstract}
We construct an explicit set of algebraically independent 
generators for the center of the universal enveloping algebra of the 
centralizer of a nilpotent matrix in the Lie algebra $\mathfrak{gl}_N(\C)$.
In particular, this gives a new proof of the freeness of the center,
 a result 
first proved by Panyushev, Premet and
Yakimova.
\end{abstract}

\maketitle

\section{Introduction}

Let $\lambda = (\lambda_1,\dots,\lambda_n)$ be a composition of $N$
such that either
$\lambda_1\geq\cdots\geq \lambda_n$ or
$\lambda_1\leq\cdots\leq \lambda_n$.
Let $\mathfrak{g}$ be the Lie algebra $\mathfrak{gl}_N(F)$, where $F$
is an algebraically closed field of characteristic zero.
Let $e \in \mathfrak{g}$ be the nilpotent matrix
consisting
of Jordan blocks
of sizes $\lambda_1,\dots,\lambda_n$ in order down the diagonal,
and write $\mathfrak{g}_e$ for the centralizer of $e$ in $\mathfrak{g}$.
Panyushev, Premet and Yakimova \cite{PPY} have recently proved that
$S(\mathfrak{g}_e)^{\mathfrak{g}_e}$,
the algebra of invariants for the adjoint
action of
$\mathfrak{g}_e$ on the symmetric algebra $S(\mathfrak{g}_e)$,
is a free polynomial algebra on $N$ generators.
Moreover, viewing $S(\mathfrak{g}_e)$ as a graded algebra
as usual so $\mathfrak{g}_e$ is concentrated in degree one,
they show that if $x_1,\dots,x_N$ are homogeneous generators for
$S(\mathfrak{g}_e)^{\mathfrak{g}_e}$ of degrees $d_1 \leq \cdots \leq d_N$,
then the sequence $(d_1,\dots,d_N)$ of {\em invariant degrees} is
equal to
\begin{align*}
(\overbrace{1,\dots,1^{\phantom{o\!\!\!}}}^{\lambda_{1} \: \text{1's}},\overbrace{2,\dots,2^{\phantom{o\!\!\!}}}^{\lambda_{2}\:\text{2's}},\dots,\overbrace{n,\dots,n^{\phantom{o\!\!\!}}}^{\lambda_{n}\:\text{$n$'s}})
\qquad\text{if $\lambda_1\geq\cdots\geq \lambda_n$,}\\
(\underbrace{1,\dots,1}_{\lambda_{n} \: \text{1's}},\underbrace{2,\dots,2}_{\lambda_{n-1}\:\text{2's}},\dots,\underbrace{n,\dots,n}_{\lambda_{1}\:\text{$n$'s}})
\qquad\text{if $\lambda_1\leq\cdots\leq \lambda_n$.}
\end{align*}
This is just one instance of the following
conjecture formulated in this generality by Premet: 
{\em For any semisimple Lie
algebra $\mathfrak{g}$
and any element $e \in \mathfrak{g}$ the
invariant algebra $S(\mathfrak{g}_e)^{\mathfrak{g}_e}$ is free}.
In \cite{PPY} this conjecture has already been
verified in many other situations besides the type $A$ case discussed here.

Returning to our special situation, 
let $Z(\mathfrak{g}_e)$ denote the center of the universal enveloping
algebra $U(\mathfrak{g}_e)$.
The standard filtration on $U(\mathfrak{g}_e)$ induces a filtration
on the 
subalgebra $Z(\mathfrak{g}_e)$ such that the associated
graded algebra $\gr Z(\mathfrak{g}_e)$ is canonically
identified with $S(\mathfrak{g}_e)^{\mathfrak{g}_e}$; see \cite[2.4.11]{Dix}.
We can lift the algebraically independent generators  $x_1,\dots,x_N$
from
$\gr Z(\mathfrak{g}_e)$ to $Z(\mathfrak{g}_e)$
to deduce (without resorting to Duflo's
theorem \cite[10.4.5]{Dix}) that $Z(\mathfrak{g}_e)$ is also a
free polynomial algebra.
The purpose of this note is to derive an explicit formula
for a set $z_1,\dots,z_N$ of algebraically independent
generators for $Z(\mathfrak{g}_e)$, generalizing the well known set of
generators of $Z(\mathfrak{g})$ itself  (the special case $e=0$) that arise
from the Capelli identity. We call these
the {\em elementary generators} for $Z(\mathfrak{g}_e)$.
Passing back down to the associated graded algebra,
one can easily obtain from them
an explicit set of {\em elementary invariants}
that generate
$S(\mathfrak{g}_e)^{\mathfrak{g}_e}$.

To formulate the main result precisely,
we must first introduce some notation for
elements of $\mathfrak{g}_e$.
Let $e_{i,j}$ denote the $ij$-matrix unit in $\mathfrak{g}$.
Draw a diagram with rows numbered $1,\dots,n$ from top to bottom and columns
numbered $1,2,\dots$ from left to right, consisting of
$\lambda_i$ boxes on the $i$th row in
columns $1,\dots,\lambda_i$, for each $i=1,\dots,n$.
Write the numbers $1,\dots,N$ into the boxes along rows, and use the notation
$\row(i)$ and $\col(i)$ for the row and column number of the box containing the
entry $i$.
For instance, if $\lambda = (4,3,2)$ then the diagram is
$$
\diagram{$1$&$2$&$3$&$4$\cr$5$&$6$&$7$\cr $8$&$9$\cr}
$$
and the nilpotent matrix $e$ of Jordan type $\lambda$ is
equal to $e_{1,2}+e_{2,3}+e_{3,4}+e_{5,6}+e_{6,7}+e_{8,9}$.
For $1 \leq i,j \leq n$ and $\lambda_j - \min(\lambda_i,\lambda_j)
\leq r < \lambda_j$, define
\begin{equation}\label{g_e_basis}
  e_{i,j;r} := \sum_{\substack{1 \leq h,k \leq N \\ \row(h) = i, \row(k) = j \\ \col(k)-\col(h) = r}} e_{h,k}.
\end{equation}
The vectors $\{e_{i,j;r}\mid 1 \leq i,j \leq n, \lambda_j - \min(\lambda_i,\lambda_j) \leq r < \lambda_j\}$ form a basis for $\mathfrak{g}_e$;
see \cite[Lemma 7.3]{BKshifted}.
Write $\mu \subseteq \lambda$ if $\mu = (\mu_1,\dots,\mu_n)$ is a composition
with $0 \leq \mu_i \leq \lambda_i$ for each $i=1,\dots,n$.
Also let $|\mu| := \mu_1+\cdots+\mu_n$ and $\ell(\mu)$ denote the number of
non-zero parts of $\mu$.
Recall that $(d_1,\dots,d_N)$ are the invariant degrees
as defined above.
Given $0 \neq \mu \subseteq \lambda$
such that $\ell(\mu) = d_{|\mu|}$,
suppose that the non-zero parts of $\mu$ are in the entries indexed by
$1 \leq i_1 < \cdots < i_{d} \leq n$.
Define the $\mu$th {\em column determinant} 
\begin{equation}\label{away}
  \cdet(\mu) := 
    \sum_{w \in S_d} \sgn(w) \tilde e_{i_{w1}, i_1;\mu_{i_1}-1}
    \tilde e_{i_{w2}, i_2;\mu_{i_2}-1}
    \cdots \tilde e_{i_{w d}, i_{d}; \mu_{i_{d}}-1},
\end{equation}
where 
$\tilde e_{i,j;r} := e_{i,j;r}- \delta_{r,0} \delta_{i,j} (i-1) \lambda_i$.
We note by Lemma~\ref{last} below that
all the $e_{i,j;r}$'s appearing
on the right hand side of (\ref{away})
necessarily satisfy the inequality
$\lambda_j - \min(\lambda_i,\lambda_j) \leq r < \lambda_j$, so
$\cdet(\mu)$ is well-defined.
For $r=1,\dots,N$, define
\begin{equation}\label{Main_formula}
  z_r := 
  \sum_{\substack{\mu \subseteq \lambda \\ |\mu| = r, \ell(\mu) = d_r }} 
  \cdet(\mu).
\end{equation}

\vspace{2mm}
{\noindent \bf Main Theorem.} {\em
The elements $z_1,\dots,z_N$ 
are algebraically independent generators
for $Z(\mathfrak{g}_e)$.
}
\vspace{2mm}

In the situation that $\lambda_1=\cdots=\lambda_n$,
our Main Theorem was proved already by
Molev \cite{M}, following Rais and Tauvel \cite{RT} who
established the freeness
of $S(\mathfrak{g}_e)^{\mathfrak{g}_e}$ in that case.
Our proof for general $\lambda$ follows the same strategy as Molev's proof,
but we need to replace the truncated Yangians with
their shifted analogs from \cite{BKshifted}.
We have also included a self-contained proof
of the freeness of $S(\mathfrak{g}_e)^{\mathfrak{g}_e}$,
although as we have said this was already established in \cite{PPY}.
Our approach to that is similar to the argument in \cite{RT}
and different from \cite{PPY}.

One final comment.
In this introduction we have formulated the Main Theorem 
assuming either that $\lambda_1 \geq \cdots \geq \lambda_n$
or that $\lambda_1 \leq \cdots \leq \lambda_n$.
Presumably most readers will prefer the former choice.
However
in the remainder of the article we will only actually prove the results in
the latter situation, 
since that is the convention adopted in \cite{BKshifted}--\cite{BKschur}.
This is justified because the two formulations of the 
Main Theorem
are simply equivalent, by an elementary 
argument involving twisting with an antiautomorphism
of $U(\mathfrak{g})$ of the form
$e_{i,j} \mapsto -e_{i',j'} + \delta_{i,j} c$.

The remainder of the article is organized as follows.
In $\S$2, we derive a new ``quantum determinant''
formula for the central elements of
the shifted Yangians. In $\S$3 we descend from there to the
universal enveloping algebra $U(\mathfrak{g}_e)$ to prove that the elements
$z_r$ are indeed central. Then in $\S$4
we prove the freeness of $S(\mathfrak{g}_e)^{\mathfrak{g}_e}$ by restricting to a carefully chosen slice.

\section{Shifted quantum determinants}

The shifted Yangian $Y_n(\sigma)$ is defined in \cite{BKshifted}.
Here are some of the details.
Let $\sigma = (s_{i,j})_{1 \leq i,j \leq n}$ be 
an $n \times n$ {\em shift matrix}, that is, 
all its entries are non-negative integers and
$s_{i,j} + s_{j,k} = s_{i,k}$ whenever $|i-j| + |j-k| = |i-k|$.
Then $Y_n(\sigma)$ is the associative algebra over $F$ defined by generators
\begin{align*}
  &\{D_i^{(r)} \mid 1 \le i \le n, r > 0 \}, \\
  &\{E_i^{(r)} \mid 1 \le i \le n, r > s_{i,i+1}\}, \\
  &\{F_i^{(r)} \mid 1 \le i \le n, r > s_{i+1,i}\} 
\end{align*}
subject to certain relations.
See \cite[\S 2]{BKshifted} for the full set.

For $1 \le i < j \le n$ and $r > s_{i,j}$, define elements 
$E_{i,j}^{(r)} \in Y_n(\sigma)$ recursively by
\begin{equation}
  E_{i, i+1}^{(r)} := E_i^{(r)}, \hspace{5mm} 
    E_{i,j}^{(r)} := [E_{i,j-1}^{(r-s_{j-1,j})}, E_{j-1}^{(s_{j-1,j}+1)}].
\end{equation}
Similarly, for $1 \le i < j \le n$ and $r > s_{j,i}$ define elements 
$F_{i,j}^{(r)} \in Y_n(\sigma)$ by
\begin{align}\label{f1}
  F_{i, i+1}^{(r)} := F_i^{(r)},\hspace{5mm}  
    F_{i,j}^{(r)} := [F_{j-1}^{(s_{j,j-1}+1)}, F_{i,j-1}^{(r - s_{j,j-1})}].
\end{align}
As in \cite[\S 2]{BKrep}, we introduce a new set of generators 
for $Y_n(\sigma)$.  
For $1 \le i < j \le n$ define the power series 
$E_{i,j}(u), F_{i,j}(u) \in Y_n(\sigma)[[u^{-1}]]$ by
\begin{equation}\label{f2}
  E_{i,j}(u) := \sum_{r > s_{i,j}} E_{i,j}^{(r)} u^{-r}, \hspace{5mm}
    F_{i,j}(u) := \sum_{r > s_{j,i}} F_{i,j}^{(r)} u^{-r},
\end{equation}
and set $E_{i,i}(u) = F_{i,i}(u) = 1$ by convention.
Also define
\begin{equation*}
  D_i(u) := \sum_{r \ge 0} D_{i}^{(r)} u ^{-r} \in  Y_n(\sigma)[[u^{-1}]],
\end{equation*}
for $1 \le i \le n$, where $D_{i}^{(0)} = 1$ by convention.
Let $D(u)$ denote the $n \times n$ diagonal matrix with $ii$-entry $D_i(u)$ for $1 \le i \le n$, let $E(u)$ denote the $n \times n$ upper triangular matrix with $ij$-entry $E_{i,j}(u)$ for $1\le i \le j \le n$, and let $F(u)$ denote the $n \times n$ lower triangular matrix with $ji$-entry 
$F_{i,j}(u)$ for $1 \le i \le j \le n$.  Consider the product
\begin{equation*}
  T(u) = F(u) D(u) E(u)
\end{equation*}
of matrices with entries in $Y_n(\sigma)[[u^{-1}]]$.  
The $ij$-entry of the matrix $T(u)$ defines a power series 
\begin{equation}\label{T_gen}
  T_{i,j}(u) = \sum_{r \ge 0} T_{i,j}^{(r)} u^{-r} := 
    \sum_{k=1}^{\min{(i,j)}} F_{k,i}(u) D_k(u) E_{k,j}(u)
\end{equation}
for some new elements $T_{i,j}^{(r)} \in Y_n(\sigma)$.  
Note that $T_{i,j}^{(0)} = \delta_{i,j}$ and $T_{i,j}^{(r)} = 0$ 
for $0 < r \le s_{i,j}$.

If the matrix $\sigma$ is the zero matrix, we denote $Y_n(\sigma)$ simply 
by $Y_n$. The algebra $Y_n$ is the Yangian 
associated to the Lie algebra $\mf{gl}_n(F)$;
see \cite[$\S$1]{MNO} for its usual definition.
In general, by \cite[Corollary 2.2]{BKshifted}, there exists an injection 
$Y_n(\sigma) \hookrightarrow Y_n$ which sends the elements 
$D_i^{(r)}, E_i^{(r)}$, and $F_i^{(r)}$ in $Y_n(\sigma)$ to the
elements with the same name in $Y_n$.  However, this injection 
usually does not 
send {\em all} the elements 
$E_{i,j}^{(r)}, F_{i,j}^{(r)}$ and $T_{i,j}^{(r)}$ 
of $Y_n(\sigma)$ to the elements with
the same name in $Y_n$. 
For the remainder of this section we will use this injection to identify
$Y_n(\sigma)$ with a subalgebra of $Y_n$.
To avoid confusion the elements 
$E_{i,j}^{(r)}, F_{i,j}^{(r)}$, and $T_{i,j}^{(r)}$ of $Y_n(\sigma)$ 
will be denoted ${^\sigma\!}E_{i,j}^{(r)}, {^\sigma\!}F_{i,j}^{(r)}$, and 
${^\sigma}T_{i,j}^{(r)}$ respectively, while $E_{i,j}^{(r)}, F_{i,j}^{(r)}$, and 
$T_{i,j}^{(r)}$ will refer to the elements of $Y_n$.
Similarly we write ${^\sigma\!}E_{i,j}(u),
{^\sigma\!}F_{i,j}(u)$ and ${^\sigma}T_{i,j}(u)$ for the
power series (\ref{f2})--(\ref{T_gen})
computed in $Y_n(\sigma)[[u^{-1}]]$ to distinguish them from their counterparts
in $Y_n[[u^{-1}]]$.

For an $n \times n$ matrix $A = (a_{i,j})_{1 \leq i,j \leq n}$ with
entries in some associative (but not necessarily commutative) algebra,
we define its {\em column determinant}
\begin{equation}\label{cdetdef}
\cdet A := \sum_{w \in S_n} \sgn(w) a_{w1,1} a_{w2,2} \cdots a_{wn,n}.
\end{equation}
For $1 \leq j \le n$, 
we  define a {\em left $j$-minor} of $A$ to be a $j \times j$ submatrix
of the form
  \[
    \left(
    \begin{array}{cccc}
      a_{i_1, 1} & a_{i_1, 2} & \cdots & a_{i_1, j} \\
      a_{i_2, 1} & a_{i_2, 2} & \cdots & a_{i_2, j} \\
      \vdots&\vdots&\ddots&\vdots\\
      a_{i_j, 1} & a_{i_j, 2} & \cdots & a_{i_j, j} \\
    \end{array}
    \right)
  \]
  for $1\le i_1 < i_2 < \dots < i_j \le n$. The following lemma is 
an easy exercise.

\begin{Lemma}\label{L1}
If for a fixed $1 \leq j \le n$, the $\cdet$ of every left $j$-minor of an 
$n \times n$ matrix $A$ with entries in some associative algebra
is zero then $\cdet A=0$.
\end{Lemma}

By \cite[Theorem 2.10]{MNO}, it is known that
the coefficients of the power series
\begin{align} \label{cj}
  C_n(u) &:= \cdet\left(
  \begin{array}{cccc}
    T_{1,1}(u) & T_{1,2}(u-1)&\cdots&T_{1,n}(u-n+1)\\
    T_{2,1}(u) & T_{2,2}(u-1)&\cdots&T_{2,n}(u-n+1)\\
    \vdots&\vdots&\ddots&\vdots\\
    T_{n,1}(u) & T_{n,2}(u-1)&\cdots&T_{n,n}(u-n+1)
  \end{array}
  \right)\\\intertext{belong to the center
of $Y_n$. Define}
\label{C_sig1}
  {^\sigma}C_{n}(u) &:= \cdet\left(
  \begin{array}{cccc}
    {^\sigma}T_{1,1}(u) & {^\sigma}T_{1,2}(u-1)&\cdots&{^\sigma}T_{1,n}(u-n+1)\\
    {^\sigma}T_{2,1}(u) & {^\sigma}T_{2,2}(u-1)&\cdots&{^\sigma}T_{2,n}(u-n+1)\\
    \vdots&\vdots&\ddots&\vdots\\
    {^\sigma}T_{n,1}(u) & {^\sigma}T_{n,2}(u-1)&\cdots&{^\sigma}T_{n,n}(u-n+1)
  \end{array}
\right).
\end{align}
The goal in the remainder of the section
is to prove the following theorem.
Note this result is false without the assumption that $\sigma$ is upper triangular.

\begin{Theorem}\label{ke}
Assuming that the shift matrix $\sigma$ is upper triangular, i.e.
$s_{i,j} = 0$ for $i > j$, we have that 
${^\sigma}C_{n}(u) = C_{n}(u)$, equality in $Y_n[[u^{-1}]]$. 
In particular, the coefficients of the power series
${^\sigma}C_n(u)$ belong to the center of $Y_n(\sigma)$.
\end{Theorem}

For the proof, assume from now on that $\sigma$ is upper triangular.
For $0 \leq j \leq n$, let $X_j$ be the $n \times n$ matrix
whose first $j$ columns are the same as the first $j$ columns
of the matrix in (\ref{cj}) and whose last $(n-j)$ columns
are the same as the last $(n-j)$ columns of the matrix in (\ref{C_sig1}).
In this notation, the theorem asserts that
$\cdet X_0 = \cdet X_n$.
So we just need to check for each $j=1,\dots,n$ that 
\begin{equation}\label{want}
\cdet X_{j-1} = \cdet X_{j}.
\end{equation}
To see this, fix $j$ and let $v := u-j+1$ for short.
Given a column vector $\vec{a}$ of height $n$,
let $X(\vec{a})$ be the matrix obtained from $X_j$ by replacing the
$j$th column by $\vec{a}$.
For $1 \leq k \leq j$, introduce the following column vectors:
$$
\vec{a} := \left(
\begin{array}{c}
{^\sigma}T_{1,j}(v)\\
{^\sigma}T_{2,j}(v)\\
\vdots\\
{^\sigma}T_{n,j}(v)
\end{array}
\right),\quad
\vec{b}_k := \left(
\begin{array}{c}
T_{1,k}(v)\\
T_{2,k}(v)\\
\vdots\\
T_{n,k}(v)
\end{array}
\right),\quad
\vec{c}_k := \left(
\begin{array}{c}
0 \\
\vdots \\
0 \\
D_k(v)\\
F_{k,k+1}(v)D_k(v)\\
\vdots\\
F_{k,n}(v)D_k(v)\\
\end{array}
\right).
$$
Also define
$$
\left(
\begin{array}{c}
d_{1,k}\\
d_{2,k}\\
\vdots\\
d_{k-1,k}
\end{array}
\right)
:=
\left(
\begin{array}{cccc}
T_{1,1}(v)&T_{1,2}(v)&\cdots&T_{1,k-1}(v)\\
T_{2,1}(v)&T_{2,2}(v)&\cdots&T_{2,k-1}(v)\\
\vdots&\vdots&\ddots&\vdots\\
T_{k-1,1}(v)&T_{k-1,2}(v)&\cdots&T_{k-1,k-1}(v)
\end{array}
\right)^{\!\!\!-1}
\left(
\begin{array}{c}
T_{1,k}(v)\\
T_{2,k}(v)\\
\vdots\\
T_{k-1,k}(v)
\end{array}
\right)
$$
and set $e_k := {^\sigma\!}E_{k,j}(v)$.
In particular, $e_j = 1$.

\begin{Lemma}\label{l1}
$\displaystyle \vec{a} = \sum_{k=1}^j \vec{c}_k e_k$.
\end{Lemma}

\begin{proof}
In view of the assumption that $\sigma$ is upper triangular, 
we have by (\ref{f1})--(\ref{f2}) that ${^\sigma\!}F_{i,j}(v) = F_{i,j}(v)$
for all $1 \leq i\leq j \leq n$.
Now the lemma follows from the definition (\ref{T_gen}).
\end{proof}

\begin{Lemma}\label{l2}
For $1 \leq k \leq j$, we have that
$\displaystyle \vec{c}_k = \vec{b}_k - \sum_{l=1}^{k-1} \vec{b}_l d_{l,k}$.
\end{Lemma}

\begin{proof}
Take $1 \leq i \leq n$ and consider the $i$th entry of the column vectors
on either side of the equation.
If $i \geq k$ then we need to show that
$F_{k,i}(v) D_k(v) = T_{i,k}(v) - \sum_{l=1}^{k-1} T_{i,l}(v) d_{l,k}$,
which is immediate from the identity \cite[(5.4)]{BKparabpres}.
If $i < k$ then we need to show that
$0 = T_{i,k}(v) - 
\sum_{l=1}^{k-1} T_{i,l}(v) d_{l,k}$.
To see this, note by the definition of $d_{l,k}$ that
$\sum_{l=1}^{k-1} T_{i,l}(v) d_{l,k}$ is equal to the matrix product
$$
\big(
T_{i,1}(v)\: \cdots \:T_{i,k-1}(v)
\big)
\!\left(
\!\begin{array}{cccc}
T_{1,1}(v)&T_{1,2}(v)&\cdots&T_{1,k-1}(v)\\
T_{2,1}(v)&T_{2,2}(v)&\cdots&T_{2,k-1}(v)\\
\vdots&\vdots&\ddots&\vdots\\
T_{k-1,1}(v)&T_{k-1,2}(v)&\cdots&T_{k-1,k-1}(v)
\end{array}
\!\right)^{\!\!\!-1}
\!\!\!\left(
\!\begin{array}{c}
T_{1,k}(v)\\
T_{2,k}(v)\\
\vdots\\
T_{k-1,k}(v)
\end{array}
\!\right)\!.
$$
The left hand row vector is the $i$th row of the matrix being inverted,
so this product does indeed equal $T_{i,k}(v)$.
\end{proof}

\begin{Lemma}\label{l3}
For any $1 \leq k \leq j-1$ and any $f$, we have that
$\cdet X(\vec{b}_k f) = 0$.
\end{Lemma}

\begin{proof}
We apply Lemma~\ref{L1}.
Take $1 \leq i_1 < \cdots < i_j \leq n$.
The corresponding 
left $j$-minor of $X(\vec{b}_k f)$ is equal to
$$
 \left(
  \begin{array}{ccccc}
    T_{i_1,1}(u) & T_{i_1,2}(u-1) & \cdots & T_{i_1,j-1}(u-j+2) & T_{i_1,k}(u-j+1)f \\
    T_{i_2,1}(u) & T_{i_2,2}(u-1) & \cdots & T_{i_2,j-1}(u-j+2) & T_{i_2,k}(u-j+1)f \\
    \vdots & \vdots & \ddots & \vdots & \vdots \\
    T_{i_j,1}(u) & T_{i_j,2}(u-1) & \cdots & T_{i_j,j-1}(u-j+2) & T_{i_j,k}(u-j+1)f
  \end{array}
  \right).
$$
The $\cdet$ of this matrix is zero by \cite[(8.4)]{BKparabpres}.
\end{proof}

Now we can complete the proof of Theorem~\ref{ke}.
Since $\cdet$ is linear in each column, we have by Lemmas~\ref{l1}-\ref{l3} 
that
\begin{align*}
\cdet X(\vec{a}) &= 
\sum_{k=1}^j \cdet X(\vec{b}_k e_k) - \sum_{k=1}^j \sum_{l=1}^{k-1}
\cdet X(\vec{b}_l d_{l,k} e_k)= \cdet X(\vec{b}_j).
\end{align*}
Since $X_{j-1} = X(\vec{a})$ and $X_j = X(\vec{b}_j)$,
this verifies (\ref{want}) hence the theorem.

\section{The central elements $z_r$}

For the remainder of the article,
$\lambda = (\lambda_1,\dots,\lambda_n)$ denotes a fixed composition of $N$
such that $\lambda_1 \leq \cdots \leq \lambda_n$
and $\sigma = (s_{i,j})_{1 \leq i,j \leq n}$ denotes the 
upper triangular shift matrix
defined by
$s_{i,j} := \lambda_j - \min(\lambda_i,\lambda_j)$.
Let $\mathfrak{g} = \mathfrak{gl}_N(F)$ and $e \in \mathfrak{g}$ be the 
nilpotent matrix consisting Jordan blocks 
of sizes $\lambda_1,\dots,\lambda_n$ down the diagonal.
Recall from the introduction 
that the centralizer $\mathfrak{g}_e$ of $e$ in $\mathfrak{g}$
has basis
\begin{equation}\label{thebasis}
\{e_{i,j;r}\mid 1 \leq i,j \leq n, s_{i,j} \leq r
< \lambda_j\}
\end{equation}
where $e_{i,j;r}$ is the element defined by (\ref{g_e_basis}).
We will view $\mathfrak{g}_e$ as a $\Z$-graded Lie algebra by declaring
that the basis element $e_{i,j;r}$ is of degree $r$.
There is an induced $\Z$-grading on the universal enveloping algebra
$U(\mathfrak{g}_e)$.

In this section we are going to prove that the
elements
$z_1,\dots,z_N$ of $U(\mathfrak{g}_e)$ 
from (\ref{Main_formula}) 
actually belong to
the center $Z(\mathfrak{g}_e)$ of $U(\mathfrak{g}_e)$ 
by exploiting the
relationship between $U(\mathfrak{g}_e)$ and the
{\em finite $W$-algebra} $W(\lambda)$
associated to $e$.
According to the definition followed here,
$W(\lambda)$ is the quotient of the
shifted Yangian $Y_n(\sigma)$ by the two-sided ideal generated by the
elements $\left\{ D_1^{(r)}\:\big|\: r > \lambda_1 \right\}$.
This is not the usual
definition of the finite $W$-algebra, but it is equivalent to the usual
definition thanks to the main result of 
\cite{BKshifted}.
The notation $T_{i,j}^{(r)}$ will from now on denote the canonical
image 
in the quotient algebra $W(\lambda)$
of the element $T_{i,j}^{(r)} \in Y_n(\sigma)$ from (\ref{T_gen})
(which was also denoted ${^\sigma}T_{i,j}^{(r)}$ in the previous section).
Recall that $T_{i,j}^{(r)} = 0$ for $0 < r < s_{i,j}$.
In addition, now that we have passed to the quotient $W(\lambda)$,
the following holds by
\cite[Theorem 3.5]{BKrep}.

\begin{Lemma}
$T_{i,j}^{(r)} = 0$
for all 
$r > \lambda_j$.
\end{Lemma}

So the power series $T_{i,j}(u) := \sum_{r \geq 0} T_{i,j}^{(r)} u^{-r}
\in W(\lambda)[[u^{-1}]]$
is actually a polynomial and
$u^{\lambda_j} T_{i,j}(u)$ belongs to $W(\lambda)[u]$.
Hence
\begin{equation*}
\cdet
\left(
\!  \begin{array}{cccc}
    \!u^{\lambda_1}T_{1,1}(u) \!& \!(u-1)^{\lambda_2}T_{1,2}(u-1)\!&\cdots&\!(u-n+1)^{\lambda_n}T_{1,n}(u-n+1)\!\\
    \!u^{\lambda_1}T_{2,1}(u) \!& \!(u-1)^{\lambda_2}T_{2,2}(u-1)\!&\cdots&\!(u-n+1)^{\lambda_n}T_{2,n}(u-n+1)\!\\
    \vdots&\vdots&\ddots&\vdots\\
    \!u^{\lambda_1}T_{n,1}(u) \!& \!(u-1)^{\lambda_2}T_{n,2}(u-1)\!&\cdots&\!(u-n+1)^{\lambda_n}T_{n,n}(u-n+1)\!
\end{array}\!\right)
\end{equation*}
gives us a well-defined polynomial $Z(u) \in W(\lambda)[u]$.
We have that
\begin{equation}\label{zu}
Z(u) = u^N + Z_1 u^{N-1} + \cdots + Z_{N-1} u + Z_N
\end{equation}
for elements $Z_1,\dots,Z_N \in W(\lambda)$.

\begin{Lemma}\label{midterm}
The elements $Z_1,\dots,Z_N$
belong to the center $Z(W(\lambda))$ of $W(\lambda)$.
\end{Lemma}

\begin{proof}
This follows from Theorem~\ref{ke},
because $Z(u)$ is equal to the canonical image of
the power series from (\ref{C_sig1})
multiplied by
$u^{\lambda_1} (u-1)^{\lambda_2} \cdots (u-n+1)^{\lambda_n}$.
\end{proof}
 
We define a filtration 
$\F_0 W(\lambda) \subseteq
\F_1 W(\lambda) \subseteq \cdots$
on $W(\lambda)$, which we call the {\em loop filtration},
by declaring that each generator $T_{i,j}^{(r+1)}$ is of filtered 
degree $r$.
In other words, $\F_r W(\lambda)$ is the span of all monomials
of the form $T_{i_1,j_1}^{(r_1+1)} \cdots T_{i_k, j_k}^{(r_k+1)}$
such that $r_1+\cdots+r_k \leq r$.
For an element $x\in \F_r W(\lambda)$, we write
$\gr_r x$ for the canonical image of $x$ in the $r$th graded
component of the associated graded algebra 
$\gr W(\lambda)$.
Applying the PBW theorem for $W(\lambda)$ from \cite[Lemma 3.6]{BKrep}, it 
follows that the loop filtration as defined here coincides with 
the filtration defined at the beginning of \cite[$\S$3]{BKschur}.
So we can restate \cite[Lemma 3.1]{BKschur} as follows.

\begin{Lemma}\label{c}
There is a unique isomorphism of graded algebras
$$
\varphi:\gr W(\lambda) \stackrel{\sim}{\longrightarrow} U(\mathfrak{g}_e)
$$
such that
$\varphi(\gr_{r} T_{i,j}^{(r+1)})= (-1)^r e_{i,j;r}$
for all $1 \leq i,j \leq n$ and $s_{i,j} \leq r < \lambda_j$.
\end{Lemma}

Let $(d_1,\dots,d_N)$ be the sequence of invariant degrees
defined in the first paragraph of the introduction.
Recall also the elements $z_1,\dots,z_N$ of $U(\mathfrak{g}_e)$
defined by the equation (\ref{Main_formula}).
The goal in the remainder of the section is to prove the following theorem.

\begin{Theorem}\label{gt}
For $r=1,\dots,N$, the element $Z_r \in Z(W(\lambda))$ 
belongs to $\F_{r-d_r} W(\lambda)$ and
$\varphi(\gr_{r-d_r} Z_r) = (-1)^{r-d_r} z_r$.
In particular, the elements $z_1,\dots,z_N$ belong to the center
$Z(\mathfrak{g}_e)$ of $U(\mathfrak{g}_e)$.
\end{Theorem}

To prove the theorem, we begin with several lemmas.

\begin{Lemma}\label{g1}
For $r =1,\dots,N$, we have that
$$
Z_r = 
\sum_{\substack{\mu \subseteq \lambda \\ |\mu| = r}}
\sum_{\nu \subseteq \mu}
\left[\left(
\prod_{i=1}^n (1-i)^{\mu_i-\nu_i}
\binom{\lambda_i-\nu_i}{\lambda_i-\mu_i}
\right)\times \left(\sum_{w \in S_n}
\sgn(w) 
T_{w1,1}^{(\nu_1)} \cdots T_{wn,n}^{(\nu_n)}
\right)\right].
$$
\end{Lemma}

\begin{proof}
Before we begin, we point out that when $i=1$ the term
$(1-i)^{\mu_i-\nu_i}$ in the product on the right hand side should
be interpreted as $1$ if $\nu_1=\mu_1$ and as $0$ otherwise.
Write $\coeff_r \left(f(u)\right)$ for the $u^r$-coefficient of a polynomial
$f(u)$.
By the definitions (\ref{cdetdef}) and (\ref{zu}), we have that
\begin{align*}
Z_r &= \sum_{w\in S_n}\! \sgn(w) 
\,\coeff_{N-r} \left(
u^{\lambda_1} T_{w1,1}(u) 
\times \cdots \times (u-n+1)^{\lambda_n}
T_{wn,n}(u-n+1)
 \right)\\
&=
\sum_{\substack{\mu \subseteq \lambda \\ |\mu|=r}}
\!\sum_{w \in S_n} \!\sgn(w)
\,\coeff_{\lambda_1-\mu_1}(u^{\lambda_1} T_{w1,1}(u))
\times \cdots\times
\end{align*}
\vspace{-8mm}
$$
\qquad\qquad\qquad\qquad\qquad\qquad\qquad
\coeff_{\lambda_n-\mu_n}((u-n+1)^{\lambda_n} T_{wn,n}(u-n+1)).
$$
Moreover for $i=1,\dots,n$  we have that
$$
\coeff_{\lambda_i-\mu_i}((u-i+1)^{\lambda_i} T_{wi,i}(u-i+1))
=
\sum_{\nu_i=0}^{\mu_i}
(1-i)^{\mu_i-\nu_i}
\binom{\lambda_i-\nu_i}{\lambda_i-\mu_i}
T_{wi,i}^{(\nu_i)}.
$$
Substituting into the preceding formula for $Z_r$ gives the conclusion.
\end{proof}

\begin{Lemma}\label{g3}
Suppose $\mu = (\mu_1,\dots,\mu_n)$ and
$\nu = (\nu_1,\dots,\nu_n)$ are compositions
with $\nu \subseteq \mu$.
We have that $|\nu|-\ell(\nu) \leq |\mu|-\ell(\mu)$
with equality if and only if
for each $i=1,\dots,n$ either
$\nu_i = \mu_i$ or $\nu_i =0=\mu_i-1$.
\end{Lemma}

\begin{proof}
Obvious.
\end{proof}

\begin{Lemma}\label{g2}
For $r=1,\dots,N$, we have that 
$d_r = \min \{\ell(\mu)\mid \mu\subseteq \lambda, |\mu| = r\}$.
\end{Lemma}

\begin{proof}
Set $d := d_r$ and $s := r - \lambda_n - \lambda_{n-1}-\cdots-\lambda_{n-d+2}$.
By the definition of $d_r$, we have that $1 \leq s \leq \lambda_{n-d+1}$.
The sum of the $(d-1)$ largest parts of $\lambda$ is $\lambda_n+\lambda_{n-1}+\cdots+\lambda_{n-d+2}$, which is $< r$.
Hence we cannot find $\mu \subseteq \lambda$ with $|\mu|=r$ and
$\ell(\mu) < d$.
On the other hand, $\mu := (0,\dots,0,s,\lambda_{n-d+2},\dots,\lambda_{n-1},\lambda_n)$ is a composition with $\mu \subseteq \lambda$ with $|\mu| = r$
and $\ell(\mu) = d$.
\end{proof}

\begin{Lemma}\label{last}
Given  $0 \neq \mu \subseteq \lambda$ with
$\ell(\mu) = d_{|\mu|}$, let $1 \leq i_1 < \cdots < i_d \leq n$
index the non-zero parts of $\mu$.
Then for any $w \in S_d$ and $j=1,\dots,d$ 
we have that
$\mu_{i_j} > \lambda_{i_j} - \min(\lambda_{i_{wj}},\lambda_{i_j})$.
\end{Lemma}

\begin{proof}
If ${wj} \geq j$, this is clear since the right hand side of the inequality is zero.
So suppose that $wj < j$, when the right hand side of the inequality
equals $\lambda_{i_j} - \lambda_{i_{wj}}$.
Assume for a contradiction that
$\mu_{i_j} \leq \lambda_{i_j} - \lambda_{i_{w j}}$.
Then we have that
\[
  | \mu | = \sum_{k=1}^d \mu_{i_k} \le \left(\sum_{k=1}^d \lambda_{i_k} \right)
   - \lambda_{i_{w j}}.
\]
Since $i_{w j} = i_k$ for some $k=1,\dots,d$, this implies that there exists
a composition $\nu\subseteq \lambda$ with $|\nu| = |\mu|$
and $\ell(\nu) = d-1$.
This contradicts Lemma~\ref{g2}.
\end{proof}

Now we can prove the theorem.
The term $T_{w1,1}^{(\nu_1)} \cdots T_{wn,n}^{(\nu_n)}$
in the expansion of $Z_r$ from Lemma~\ref{g1}
belongs to $\F_{|\nu|-\ell(\nu)} W(\lambda)$.
If $\nu \subseteq \mu \subseteq \lambda$ and $|\mu|=r$, 
Lemmas~\ref{g3}--\ref{g2} imply that
$|\nu|-\ell(\nu) \leq |\mu| - \ell(\mu) = r - \ell(\mu)
\leq r-d_r$.
This shows that $Z_r$ belongs to $\F_{r-d_r} W(\lambda)$.
Moreover, to compute $\gr_{r-d_r} Z_r$ we just need to 
consider the terms in the expansion of $Z_r$ that have
$\ell(\mu) = d_r$ and for each $i=1,\dots,n$ either $\nu_i=\mu_i$ or
$\nu_i=0=\mu_i-1$.
We deduce from Lemma~\ref{g1} that
$$
\gr_{r-d_r} Z_r = 
\sum_{\substack{\mu \subseteq \lambda \\ |\mu|=r,\ell(\mu) = d_r}}
\sum_{w \in S_n} \sgn(w) 
\gr_{r-d_r}
\left(
\widetilde{T}_{w1,1}^{(\mu_1)} \cdots
\widetilde{T}_{wn,n}^{(\mu_n)}
\right)
$$
where
$\widetilde{T}_{i,j}^{(r)} := T_{i,j}^{(r)} + \delta_{i,j}\delta_{r,1} (1-i)\lambda_i$.
Since $\widetilde{T}_{wi,i}^{(0)} = 0$ unless $wi=i$, we can further simplify 
this expression as follows. 
Let $d := d_r$ for short and given $\mu\subseteq \lambda$ with  $\ell(\mu) = d$
define $1 \leq i_1 < \cdots < i_d \leq n$ so that
$\mu_{i_1},\dots,\mu_{i_d} \neq 0$. Then
$$
\gr_{r-d} Z_r = 
\sum_{\substack{\mu \subseteq \lambda \\ |\mu|=r,\ell(\mu) = d}}
\sum_{w \in S_d} \sgn(w)
\left(\gr_{\mu_{i_1}-1} \widetilde{T}_{i_{w1},i_1}^{(\mu_{i_1})}\right) \cdots
\left(\gr_{\mu_{i_d}-1} \widetilde{T}_{i_{wd},i_d}^{(\mu_{i_d})}\right).
$$
Finally applying the isomorphism $\varphi$ from Lemma~\ref{c}, 
we get that
$$
\varphi(\gr_{r-d} Z_r) = \!\!\!\sum_{\substack{\mu\subseteq\lambda \\ 
|\mu|=r, \ell(\mu)=d}}
\!\!\!\sum_{w \in S_d}
\sgn(w)
(-1)^{\mu_{i_1}-1}
\tilde e_{i_{w1},i_1;\mu_{i_1}-1}
\cdots
(-1)^{\mu_{i_d}-1}
\tilde e_{i_{wd},i_d;\mu_{i_d}-1}
$$
where $\tilde e_{i,j;r} := \gr_r \widetilde T_{i,j}^{(r+1)}$.
The right hand side is $(-1)^{r-d} z_r$ according to the definitions
in the introduction. Noting finally that, since 
$Z_r$ is central in $W(\lambda)$ by Lemma~\ref{midterm}, the element
$\gr_{r-d} Z_r$ is central in $\gr W(\lambda)$,
this completes the proof of Theorem~\ref{gt}.

\section{Proof of the main theorem}

Now we consider the standard filtration on
the universal enveloping algebra $U(\mathfrak{g}_e)$
and the induced filtration on the subalgebra $Z(\mathfrak{g}_e)$.
By the PBW theorem, the associated graded algebra
$\gr U(\mathfrak{g}_e)$ is identified with the 
symmetric algebra $S(\mathfrak{g}_e)$ 
(generated by $\mathfrak{g}_e$ in degree one).
For $r=1,\dots,N$, it is immediate from the definition (\ref{Main_formula})
that the central element $z_r \in U(\mathfrak{g}_e)$
is of filtered degree $d_r$.
Let $x_r := \gr_{d_r} z_r \in S(\mathfrak{g}_e)^{\mathfrak{g}_e}$.
Explicitly,
\begin{equation}\label{xexp}
  x_r 
=  \sum_{\substack{\mu \subseteq \lambda \\ |\mu| = r,\: \ell(\mu) = d_r }} 
  \sum_{w \in S_d} \sgn(w) e_{i_{w1}, i_1;\mu_{i_1}-1}
  \cdots e_{i_{w d}, i_{d}; \mu_{i_{d}}-1} \in S(\mathfrak{g}_e)
\end{equation}
where as usual 
$1 \le i_1 < i_2 \dots < i_{d} \leq n$ denote the positions of the non-zero
entries of $\mu$.

\begin{Theorem}\label{ic}
The elements $x_1,\dots,x_N$ are algebraically independent generators
for $S(\mathfrak{g}_e)^{\mathfrak{g}_e}$.
\end{Theorem}

For the proof,
let us from now on identify $S(\mathfrak{g}_e)$ with
$F[\mathfrak{g}_e^*]$, the coordinate algebra of the affine variety
$\mathfrak{g}_e^*$.
Let
\begin{equation}\label{thfbasis}
\{f_{i,j;r}\mid 1 \leq i,j \leq n, s_{i,j} \leq r
< \lambda_j\}
\end{equation}
be the basis for $\mathfrak{g}_e^*$ that is dual to the
basis (\ref{thebasis}).
By convention, we interpret 
$f_{i,j;r}$ as zero if $r < s_{i,j}$.
The coadjoint action $\ad^*$ 
of $\mathfrak{g}_e$ on $\mathfrak{g}_e^*$ is given explicitly
by the formula
\begin{equation}\label{adstar}
(\ad^* e_{i,j;r}) (f_{k,l;s}) =
\delta_{j,l} f_{k,i;s-r} - \delta_{i,k} f_{j,l;s-r}.
\end{equation}
The induced action of $\mathfrak{g}_e$ on
$F[\mathfrak{g}_e^*]$ is defined by
$(x \cdot \theta)(y) = -\theta((\ad^* x)(y))$ for $x \in \mathfrak{g}_e,
y \in \mathfrak{g}_e^*, \theta \in F[\mathfrak{g}_e^*]$.
It is for this action that
the invariant subalgebra $S(\mathfrak{g}_e)^{\mathfrak{g}_e}$ is identified
with $F[\mathfrak{g}_e^*]^{\mathfrak{g}_e}$.
Introduce the affine subspace 
\begin{equation}
S:= f+V
\end{equation}
of $\mathfrak{g}_e^*$,
where
$f := 
f_{1,2;\lambda_2-1}+f_{2,3;\lambda_3-1}+\cdots+f_{n-1,n;\lambda_n-1}$ 
and
$V$ is the $N$-dimensional linear subspace spanned by the vectors
$\{f_{n,i;r}\mid 1 \leq i \leq n, 0 \leq r < \lambda_{i}\}$.
Let
\begin{equation}
\rho:F[\mathfrak{g}_e^*]^{\mathfrak{g}_e} \rightarrow 
F[S]
\end{equation}
be the homomorphism
defined by restricting functions 
from $\mathfrak{g}_e^*$ to the slice $S$.

\begin{Lemma}\label{c1}
The elements $\rho(x_1),\dots,\rho(x_N)$ are algebraically independent
generators of
$F[S]$.
\end{Lemma}

\begin{proof}
Take an arbitrary vector
$$
v = 
f_{1,2;\lambda_2-1}+f_{2,3;\lambda_3-1}+\cdots+f_{n-1,n;\lambda_n-1} + \sum_{j=1}^n \sum_{t=0}^{\lambda_j-1} a_{j,t} f_{n,j;t} \in S.
$$
Since $S \cong \mathbb{A}^N$, the 
algebra $F[S]$ is freely generated by the coordinate functions
$p_{j,t}:v \mapsto a_{j,t}$ for $1 \leq j \leq n$ and $0 \leq t < \lambda_j-1$.
Also note for any $1 \leq i,j \leq n$ and 
$s_{i,j} \leq r <  \lambda_j$ that
\begin{equation}\label{form}
e_{i,j;r}(v) = \left\{
\begin{array}{ll}
a_{j,r} & \text{if $i=n$,}\\
1&\text{if $j=i+1$ and $r=\lambda_j-1$,}\\
0&\text{otherwise.}
\end{array}
\right.
\end{equation}
Now fix $1 \leq r \leq N$.
Let $d := d_r$ and 
$s:= r-\lambda_n-\lambda_{n-1}-\cdots-\lambda_{n-d+2}$,
so that $1 \leq s \leq \lambda_{n-d+1}$.
We claim that $x_r(v) = (-1)^{d-1} a_{n-d+1,s-1}$,
hence $\rho(x_r) = (-1)^{d-1} p_{n-d+1,s-1}$.
Since every coordinate function $p_{j,t}$ arises in this way
for a unique $r$,
the lemma clearly follows from this claim.

To prove the claim, 
suppose we are given  $w \in S_d$
and $\mu\subseteq \lambda$
such that $|\mu|=r$,
the
non-zero parts of $\mu$ are in positions $1 \leq i_1 < \cdots < i_d \leq n$, 
and the 
monomial $e_{i_{w1},i_1;\mu_{i_1}-1}\cdots e_{i_{wd},i_d;\mu_{i_d}-1}$ from the right hand side of
(\ref{xexp}) is non-zero on $v$.
For at least one $j=1,\dots,d$, we must have that
$wj \geq j$, and for such a $j$
the fact that $e_{i_{wj}, i_j; \mu_{i_j}-1}(v) \neq 0$ implies by (\ref{form})
that
$wj=d$ and $i_d = n$.
For all other $j \neq k \in \{1,\dots,d\}$, 
we have that $wk \neq d$ hence $i_{wk} \neq n$.
But then the fact that $e_{i_{wk}, i_k; \mu_{i_k}-1}(v) \neq 0$
implies by (\ref{form}) that $i_k = i_{wk} +1$ and $\mu_{i_k} = \lambda_{i_k}$.
So $wj = d$ and $wk = k-1$ for all $k \neq j$,
which means that $w = (d\:\:d\!-\!1\:\cdots\:1)$ and $j=1$.
Moreover, $i_2=i_1+1,i_3=i_2+1,\dots,i_{d}=i_{d-1}+1=n$,
which means that $(i_1,\dots,i_d) = (n-d+1,\dots,n-1,n)$.
Hence, $\mu = (0,\dots,0,s,\lambda_{n-d+2},\dots,\lambda_{n-1},\lambda_n)$.
For this $w$ and $\mu$ it is indeed the case that
$e_{i_{w1},i_1;\mu_{i_1}-1}\cdots 
e_{i_{wd},i_1;\mu_{i_d}-1} (v) = a_{n-d+1,s-1}$ by (\ref{form}) once more.
Since $\ell(w) = d-1$ this and the definition (\ref{xexp}) implies the claim.
\end{proof}

\begin{Lemma}\label{c2}
$\rho$ is an isomorphism.
\end{Lemma}

\begin{proof}
Lemma~\ref{c1} implies that $\rho$ is surjective, so it just
remains to prove that it is injective.
Let $G := GL_N(F)$ acting naturally on $\mathfrak{g}$
by conjugation.
Let $G_e$ be the centralizer of $e$ in $G$
and identify $\mathfrak{g}_e$ with the Lie algebra of $G_e$,
i.e. tangent space $T_{\iota}(G_e)$ to $G_e$ at the identity element $\iota$,
as usual.
Considering
the coadjoint action $\Ad^*$ of $G_e$ on $\mathfrak{g}_e^*$,
we have 
that $F[\mathfrak{g}_e^*]^{G_e} = F[\mathfrak{g}_e^*]^{\mathfrak{g}_e}$.
To prove that $\rho:F[\mathfrak{g}_e^*]^{G_e}
\rightarrow F[S]$ is injective, it suffices to prove that $(\Ad^* G_e)(S)$ is
dense in $\mf{g}_e^*$,
i.e. that the map $\phi: G_e \times S \rightarrow \mf{g}_e^*,
(g,x) \mapsto (\Ad^* g)(x)$ is dominant.  
This follows if we can check that its differential $d\phi_{(\iota,f)}$
at the point $(\iota,f)$
is surjective; see e.g. \cite[Theorem 4.3.6(i)]{Springer}.
Identify the tangent spaces $T_f(S)$ and $T_f(\mathfrak{g}_e^*)$ with
$V$ and $\mathfrak{g}_e^*$.
Then the differential 
$d \phi_{(\iota,f)}:\mathfrak{g}_e \oplus V \rightarrow \mathfrak{g}_e^*$ 
is given explicitly by the map $(x,v) \mapsto (\ad^* x)(f) + v$.
We show that it is surjective by checking that
every basis element $f_{i,j;r}$ from (\ref{thfbasis}) 
belongs to its image.

To start with, it is easy to see each $f_{n,i;r}$ belongs to the
image of $d \phi_{(\iota,f)}$, since each of these vectors belongs to $V$.
Next, suppose that $1 \leq i \le j <n$ and $0 \leq r < \lambda_i$.
By (\ref{adstar}), we have that
$$
(\ad^* e_{i,j+1;\lambda_{j+1}-r-1})(f) = 
f_{j,i;r}-f_{j+1,i+1;\lambda_{i+1}-\lambda_{j+1}+r}.
$$
Using this, we get that all $f_{j,i;r}$ with $i \leq j$
belong to the image of $d \phi_{(\iota,f)}$.
Finally, suppose that $n \geq i > j \geq 1$ and $\lambda_i - \lambda_j \le r < \lambda_i$.
By (\ref{adstar}) again, we have that
\[
  (\ad^* e_{i-1,j; \lambda_i -r - 1})(f) =
  \begin{cases}
    -f_{j,i;r}, &\text{if $j=1$} \\
    f_{j-1,i-1;\lambda_j - \lambda_i + r}-f_{j,i;r}, &\text{if $j > 1$.}
  \end{cases}
\]
From this we see that all $f_{j,i;r}$'s with $i > j$ belong to the
image of $d\phi_{(\iota,f)}$ too.
This completes the proof.
\end{proof}

By Lemma~\ref{c1}, $\rho(x_1),\dots,\rho(x_N)$
are algebraically independent generators
for $F[S]$. By Lemma~\ref{c2}, $\rho$ is an isomorphism.
Hence $x_1,\dots,x_N$ are algebraically independent generators
of $F[\mathfrak{g}_e^*]^{\mathfrak{g}_e}$.
This completes the proof of Theorem~\ref{ic}.
Now we can deduce the Main Theorem from the introduction.

\begin{Corollary}\label{mc}
The elements $z_1,\dots,z_N$ are algebraically independent generators
for $Z(\mathfrak{g}_e)$.
\end{Corollary}

\begin{proof}
It is obvious that 
$\gr Z(\mathfrak{g}_e) \subseteq S(\mathfrak{g}_e)^{\mathfrak{g}_e}$.
We have observed already that $z_1,\dots,z_N \in Z(\mathfrak{g}_e)$ are of filtered degrees
$d_1,\dots,d_N$ respectively, and by
Theorem~\ref{ic} the associated graded elements
are algebraically independent generators for
$S(\mathfrak{g}_e)^{\mathfrak{g}_e}$.
By a standard filtration argument (see e.g. 
the proof of \cite[Theorem 2.13]{MNO}), 
this is enough to deduce that
$z_1,\dots,z_N$ are themselves algebraically independent generators
for $Z(\mathfrak{g}_e)$. At the same time, we have 
reproved the well-known equality
$\gr Z(\mathfrak{g}_e) = S(\mathfrak{g}_e)^{\mathfrak{g}_e}$.
\end{proof}

To conclude the article, we give one application;
see \cite[Theorem 6.10]{BKrep}, \cite[Remark 2.1]{PPY}
and the footnote to \cite[Question 5.1]{P2} for other proofs
of this result.
Recall the 
central elements $Z_1,\dots,Z_N$ of $W(\lambda)$ from Lemma~\ref{midterm}.

\begin{Corollary}
The elements
$Z_1,\dots,Z_N$ are algebraically independent
generators for the center of $W(\lambda)$.
\end{Corollary}

\begin{proof}
The loop filtration on $W(\lambda)$ induces a filtration on
$Z(W(\lambda))$. Clearly we have that
$\gr Z(W(\lambda)) \subseteq Z(\gr W(\lambda))$.
By Theorem~\ref{gt}, we know that $Z_1,\dots,Z_N \in Z(W(\lambda))$ are
of filtered degrees $1-d_1,\dots,N-d_N$ 
respectively, and by
Corollary~\ref{mc}
the associated graded elements are algebraically independent
generators for $Z(\mathfrak{g}_e)$.
Hence $Z_1,\dots,Z_N$ are algebraically independent generators
for $Z(W(\lambda))$.
At the same time, we have proved that
$\gr Z(W(\lambda)) = Z(\mathfrak{g}_e)$.
\end{proof}

\end{document}